\def \dd#1{{\bf#1}}

\def\cl#1{{\cal#1}}


\def\Max{\mathop{\rm Max}}
\def\Min{\mathop{\rm Min}}


\def\ouv#1{\smash{\mathop{#1}\limits^{\lower 1pt\hbox
{$\scriptscriptstyle\circ$}}}}

\def\hfl#1#2{\smash{\mathop{\hbox to 12mm{\rightarrowfill}}
\limits^{\scriptstyle#1}_{\scriptstyle#2}}}


\long\def\eno#1#2{\par\smallskip{\bf{#1}}{\it\ {#2}}\par\medskip}

\def\stit#1{\vskip 3mm plus 1mm minus 2mm {\bf{#1}}
		\smallskip}

\font\tir=cmbx10 at 12pt

\def\ref#1#2#3#4{{\bf #1}{\ #2}{\it ,\ #3}{,\ #4}\medskip}


\def \picture #1 by #2 (#3){\midinsert \centerline 
{\vbox to #2{\hrule width #1 heigth 0pt 
depth 0pt \null \vfill \special {picture #3}}}\endinsert }

\def\scaledpicture #1 by #2 (#3 scaled #4) {{
\dimen0 =#1 \dimen1 =$2
\divide \dimen0 by 1000 \multiply \dimen0 by #4
\divide \dimen1 by 1000 \multiply \dimen1 by #4
\picture \dimen0 by \dimen1 (#3 scaled $4)}}

\def\figure #1 #2 #3 {\midinsert \vglue 3mm 
{\vbox to #3 {\hrule width 6cm height 0cm depth 0cm \vfill
{\special {picture #1 scaled #2}}}}\vglue 2mm \endinsert}

\magnification=1200

\input psfig.sty

\overfullrule=0pt

{\centerline {\tir {Hedgehogs of Hausdorff dimension 1}}}


\bigskip
\bigskip

{\centerline {Kingshook Biswas \footnote {*} {Departimento di Matematica, Universita di Pisa,
Largo Bruno Pontecorvo 5,
56127 Pisa,
Italy. E-mail: kbiswas@dm.unipi.it} }}

\bigskip
\bigskip

{\bf Abstract.} {\it We present a construction of hedgehogs for
holomorphic maps with an indifferent fixed point. We construct,
for a family of commuting non-linearisable maps, a common hedgehog
of Hausdorff dimension 1, the minimum possible. }

\bigskip

\centerline{ {\it AMS Subject Classification: 37F50}}

\bigskip

\bigskip

\stit {1. Introduction.}

\bigskip

We consider the dynamics of a holomorphic diffeomorphism $f(z) = e^{2\pi i \alpha} z + O(z^2), \, \alpha \in
\dd{R} - \dd{Q}$, defined in a neighbourhood of the indifferent irrational fixed point 0. The map $f$ is said
to be {\it linearisable}
if there is a holomorphic change of variables $z = h(w) = w + O(w^2)$ such that
$$
h^{-1} \circ f \circ h = R_\alpha
$$
in a neighbourhood of the origin, where $R_\alpha(w) = e^{2\pi i \alpha} w$ is the rigid rotation.

\medskip

The problem of {\it linearisation}, or determining when such an $f$ is linearisable, has a long and interesting
history, including the work of H. Cremer ([Cr1], [Cr2]) in the 1920's, C.L.Siegel ([Si]) in 1942, A.D.Bruno in the
1960's, and Yoccoz ([Y]), who resolved the problem of the optimal arithmetic condition on $\alpha$ for
linearisability in 1987. In this article however we are primarily concerned with the structure of invariant sets
for the dynamics near the fixed point, rather than the issue of linearisability.

\medskip

When $f$ is linearisable, $f$ behaves like the rotation by angle $2\pi\alpha$ around 0 on the maximal domain of
linearisation called the {\it Siegel disk} of $f$, which turns out to be a domain of the form $h( \{ |w| < R \})$
where $R$ is less than or equal to the radius of convergence of $h$ around 0. In this case the compacts $h( \{ |w|
\leq r \})$  for $r < R$ are clearly completely invariant under $f$. However, even in the more general case where
$f$ is not necessarily linearisable, R.Perez-Marco found completely invariant continua for $f$ which persist:

\medskip

\eno {Theorem (Perez-Marco, [PM1]).}
{Let $f(z) = e^{2\pi i \alpha}z + O(z^2), \, \alpha \in \dd{R} - \dd{Q}$ be such that $f$ and $f^{-1}$ are
defined and univalent on a neighbourhood of the closure of a Jordan domain $\Omega \subset \dd C$ containing
$0$. Then there exists a full compact connected set $K$ contained in $\overline{\Omega}$ such that $0 \in K$,
$K \cap \partial\Omega \neq \phi$ and $f(K) = f^{-1}(K) = K$.
}

\medskip

The invariant compacts $K$ thus obtained are called {\it Siegel compacta}. A classical topological theorem of
G.D.Birkhoff ([Bi]) does in fact guarantee, for planar homeomorphisms near Lyapunov unstable points, the existence
of compacts, but which are either positive or negative invariant, not necessarily totally invariant. In the
holomorphic setting one can thus improve to obtain totally invariant compacts. Moreover, though this will not be
used in what follows, if the boundary of the Jordan domain $\Omega$ is $C^1$ smooth, then Perez-Marco has shown
that such a compact $K$ is in fact unique (see [PM5]). If $K$ is not contained in the closure of a linearisation
domain it is called a {\it hedgehog}. A {\it hedgehog} is called linearisable or non-linearisable depending on
whether it contains a linearisation domain or not. We will only deal with non-linearizable hedgehogs in
this article. Perez-Marco has shown that these have no interior (see [PM5]). The structure of such hedgehogs is
topologically complex; for example, Perez-Marco shows in [PM2] that they are not locally connected at any point
different from the fixed point. Objects such as {\it combs} (homeomorphs of the product of a Cantor set and an
interval) had been expected to be found within hedgehogs. Indeed this is stated by Perez-Marco in
[PM4], and the author uses Perez-Marco's techniques of ''tube-log
Riemann surfaces'' to construct hedgehogs containing smooth combs
in [B]. These techniques of ''tube-log Riemann surfaces'', first
introduced by Perez-Marco in [PM6] and later used again in [PM3]
and [PM4], are a powerful means of constructing a wide range of nonlinearizable
dynamics with good control on the geometry of the corresponding
hedgehogs.

\medskip

In the present article, we use these techniques to construct a hedgehog having Hausdorff
dimension 1 (the minimum possible for connected planar sets).

\eno {Main Theorem.} {There exists a Cantor set $C \subset \dd{R}$ and a family of commuting holomorphic
maps $(f_t(z) = e^{2\pi i t}z + O(z^2))_{t\in C}$, all defined and univalent on a neighbourhood of the
closed unit disk $\overline{\dd{D}}$, having a common hedgehog $K \subset \overline{\dd{D}}$, such that $K$
has Hausdorff dimension 1. The maps $f_t$ are non-linearisable for $t \in C \cap \dd{R} - \dd{Q}$.
}
\bigskip

\stit {2.  The construction.}

\bigskip

\bigskip

{\bf 2.1   General idea of the construction.}

\medskip

It is strongly advised to read this section in conjunction
with [PM3] and [PM4]. As these techniques are explained in detail
in [PM3] and especially in [PM4], and the basic elements of the
construction are the same as in [PM4], to avoid repetition the
description below is shorter and slightly informal.

\medskip

We start with the linear flow of translations $(F_t(z)=z+t)_{t\in
\dd{R}}$  on the upper half-plane $\dd{H}$ ({\bf note} : for
convenience, instead of working with dynamics and dynamical
objects such as Siegel disks, hedgehogs etc defined on the unit
disk $\dd D$, we will actually work instead with their lifts to
the upper half-plane via the universal covering $E : \dd H \to \dd
D - {0}, E(z) = e^{2\pi i z})$.

\medskip

{\bf 2.1.1)   The Riemann surface $\cl S_0$.}

\medskip

The idea is 'fold' the dynamics by introducing nonlinearities in
the dynamics, such as fixed points; these result in a shrinking of
the linearisation domain. This is accomplished by means of the
following Riemann surface with distinguished charts $\cl S_0$ (see
[PM4] section 1.2.1) :

\medskip

{\hfill {\centerline {\psfig {figure=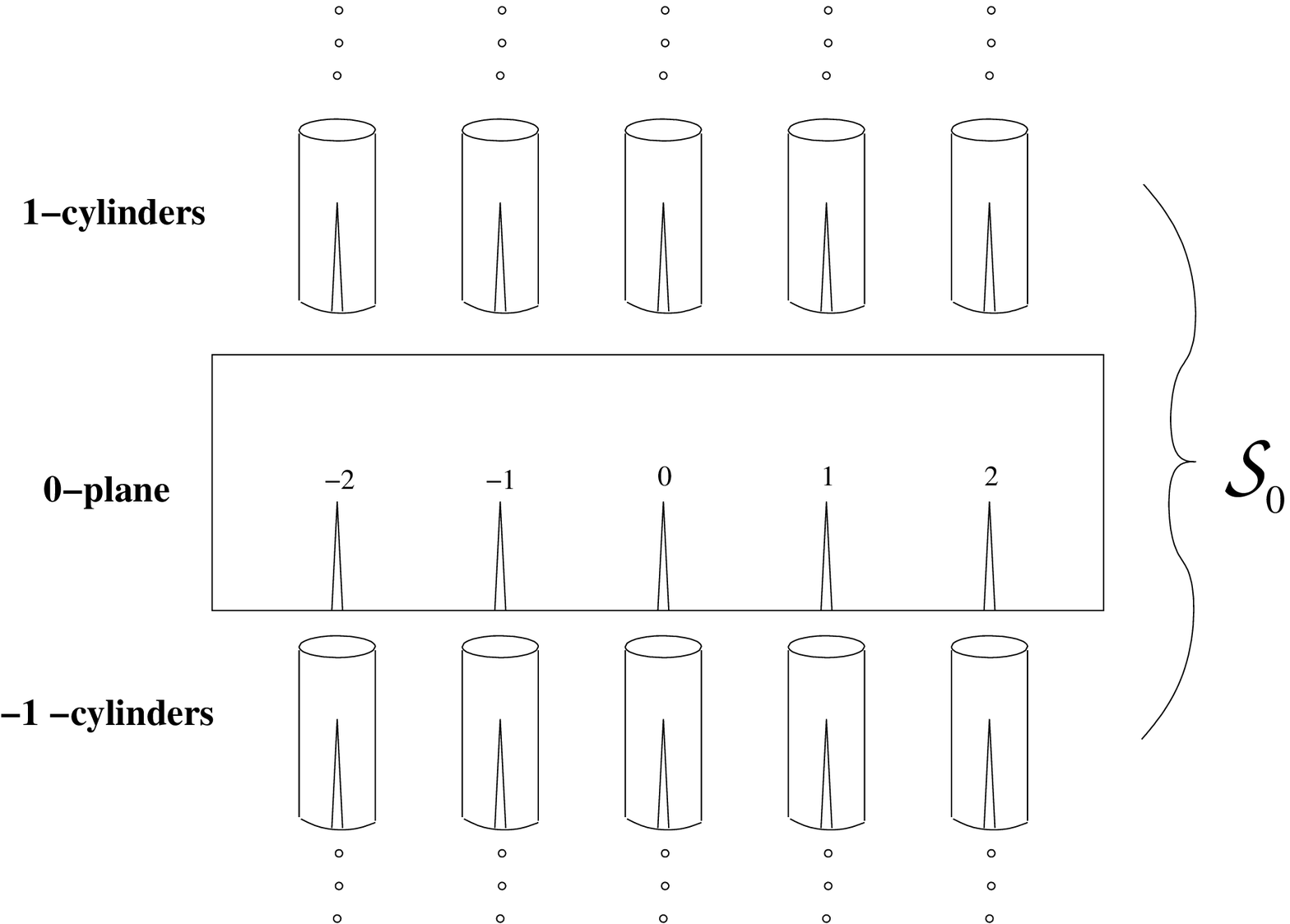,height=8cm}}}}

{\centerline {\bf Figure 1}}

\medskip

The surface $\cl S_0$ (a ''tube-log Riemann surface'') is formed from a copy of $\dd{C}$ with
'slits' $[n, n-i\infty]$ at each of the integers $n\in \dd{Z}$,
above and below each of which is pasted isometrically along the
slits a family (indexed by $\dd{Z}$) of slit cylinders
${\dd{C}}/\dd{Z}$ with slits $[0,0-i\infty]$. This construction
gives not only a Riemann surface, but also a canonical set of
charts for the surface, which allow us to write formulas for
functions defined on the surface in terms of these canonical
coordinates. We note here that the surface $\cl S_0$ has
'ramification points of infinite order' at the points
corresponding to $\dd Z$ in the $0$-plane and to $0$ in the
cylinders. These can be seen as points added in the completion of
the flat metric on $\cl S_0$ (induced from the flat metrics on its
building blocks). However it is important to note that these
points do not belong to the surface $\cl S_0$ itself, nor can the
complex structure on $\cl S_0$ be extended to these points.

\medskip

{\bf 2.1.2)   Uniformization of the Riemann surface $\cl S_0$.}

\medskip

The upper ends of the cylinders on the other hand do correspond to
points for the complex structure. Adding in these points gives a
simply connected Riemann surface which it is then not too hard to
see is biholomorphic to $\dd{C}$. So $\cl S_0$ should be
biholomorphic to the complex plane minus a doubly infinite
discrete set of points; in fact, there is an explicit formula for
the uniformization (see [PM4], section 1.2.3), which in terms of
the distinguished charts on $\cl S_0$ is given by $$ K(z) = {1
\over 2\pi i} \log ( {-\log ({1-e^{2\pi iz}})}) $$

The map $K$ satisfies the properties $$\displaylines{ K(z+1) =
K(z)+1 \quad \hbox{ on all of the 0-sheet of $\cl S_0$ } \cr
K(z)-z \to 0 \quad \hbox{ as Im } z \to + \infty \hbox{ in the
0-sheet of $\cl S_0$} \cr }$$ for the appropriate choices of the
branches of the logarithms.

\medskip

The images under $K$ of horizontal lines just above the real line
$\dd{R}$ in the 0-plane are thus 1-periodic (with respect to
parametrization by the $x$-coordinate) curves in $\dd{C}$; since
the ramification points at $\dd Z$ in the $0$-plane correspond to
ends at infinity under the uniformization, these curves oscillate
with amplitudes which tend to infinity as we take lines tending to
the real line. The upper ends of the cylinders correspond under
$K$ to a discrete set of points $P \subset \dd C$; amongst these,
those that correspond to the upper ends of the 1-cylinders lie
below, but close to, the image of the real line. 





\medskip

{\bf 2.1.3)   Lifting and folding of dynamics.}

\medskip

Consider the linear flow of translations $(F_t(z)=z+t)_{t\in
\dd{R}}$  on the upper half-plane $\dd{H}$. The upper half-plane
of the $0$-plane of $\cl S_0$ projects univalently onto the upper
half-plane $\dd{H}$, so
 we can lift each $F_t$ to a map $\tilde{F}_t(z)$ from the upper half-plane of the $0$-plane of $\cl S_0$ to
itself, given in terms of the distinguished chart by the same
formula $\tilde{F}_t(z) = z+t$. We would like to extend the lift
analytically  to a univalent map $\tilde{F}_t$ of the surface to
itself; however, it is not possible to extend to the points $\dd Z
- t$ of the $0$-plane or the points $0 - t$ of the cylinders, as
these would have to be mapped to the ramification points, and
moreover there is a non-trivial monodromy when continuing
analytically around these points.

\medskip

To overcome this problem, we can remove from the surface a small
$\epsilon$-neighbourhood (for the flat metric on $\cl S_0$) $\cl
I(\epsilon) \subset \cl S_0$ of the ramification points. On $\cl
S_0 - \cl I(\epsilon)$, for $|t| < \epsilon$, the aforementioned
problem does not arise for the maps $\tilde{F}_t$, each of which
then extends to a univalent map $\tilde{F}_t : \cl S_0 - \cl
I(\epsilon) \to \cl S_0$, which is expressed in terms of the
distinguished charts on $\cl S_0$ by the same formula
$\tilde{F}_t(z) = z+t$ (for this first step it suffices in fact to
remove only a set of horizontal slits $\dd Z + (\epsilon,
\epsilon)$ near the ramification points, but later in order to
iterate the 'folding' procedure we will have to remove
neighbourhoods of the ramification points).

\medskip

Let $T : \cl S_0 \to \cl S_0$ be the automorphism given in the
charts on $\cl S_0$ by $T(z) = z+1$. We note that the lifts of the
integer translations $z \mapsto z + n$ extend to all of $\cl S_0$
as the automorphisms $T^n : \cl S_0 \to \cl S_0$. These
automorphisms commute with the maps $\tilde{F}_t$, hence we can
also define for times $t' = n + t \in \dd Z + (-\epsilon,
\epsilon)$ the lifts $\tilde{F}_{t'} : \cl S_0 - \cl I(\epsilon)
\to \cl S_0$ by $\tilde{F}_{n + t} := T^n \circ \tilde{F}_t$.
Hence we obtain a semi-flow $(\tilde{F}_t)$ defined for times $t
\in \dd Z + (-\epsilon,\epsilon)$ on a large part $\cl S_0 - \cl
I(\epsilon)$ of the surface.

\medskip

The 'upper ends' of the cylinders on the surface are points for
the complex structure, which become fixed points (modulo the
automorphism $T : \cl S_0 \to \cl S_0$) for the dynamics on the
surface. Hence lifting the dynamics to the surface has the effect
not only of extending its domain of definition, but also of
introducing nonlinearities, such as fixed points.

\medskip

To see this dynamics on the surface as dynamics on the plane, we
use the uniformization $K : \cl S_0 \to \dd C - P$ :

\medskip

Conjugating the lifts $\tilde{F}_t$ by $K$ we obtain dynamics $K
\circ \tilde{F}_t \circ K^{-1}$ on a large part $K(\cl S_0 - \cl
I(\epsilon)) \subset \dd C$ of the plane. The dynamics flows along
the foliation shown in Figure 2, and has as fixed points (modulo
$\dd Z$) the points $P$ corresponding under $K$ to the upper ends
of the cylinders. This dynamics is what we mean by the 'folded'
dynamics.

\medskip

A key property of the construction is that small
neighbourhoods of the infinite ramification points in $\cl S_0$
correspond under $K_0$ to domains in the plane with large negative
imaginary part; any half-plane $\{ \hbox{ Im } z \geq -y \}
\subset \dd C$ where $y \geq 0$ is contained in $K(\cl S_0 - \cl
I(\epsilon))$ for $\epsilon$ small enough (see [PM4], section
1.2.4), and so we can always ensure a large domain of definition
for the dynamics.

\medskip

An equally important property of the
uniformization $K_0$ is that the lower ends of the cylinders and
points with large negative imaginary part in the $0$-plane in $\cl
S_0$ correspond to points with large negative imaginary part in
the plane. As explained later in section 2.1.5, these two
properties make it possible to iterate the 'folding' construction
indefinitely all the while maintaining a common domain of
definition, so that one can obtain dynamics in the limit.

\medskip

{\bf 2.1.4)   Normalization of the uniformization.}

\medskip

The images under $K$ of horizontal lines in the upper half-plane
of the $0$-sheet of $\cl S_0$ are graphs over the $x$-axis which
are 1-periodic with respect to parametrization by the
$x$-coordinate. For the image of a horizontal line at a height
$\delta > 0$ in the $0$-sheet, its peaks and troughs occur at the
points $K(n+1/2 + i\delta)$ and $K(n + i\delta)$ respectively
(where $n$ runs over all integers). The amplitude of oscillation
(the difference in height between peaks and troughs) is a strictly
decreasing function of $\delta$ which tends to $+\infty$ as
$\delta$ tends to $0$, and tends to $0$ as $\delta$ tends to
$+\infty$.

\medskip

We define a family of normalized uniformizations as follows:

\medskip

First we place the ramification points of $\cl S_0$ slightly below
the real axis instead of on the real axis; more precisely, for
$\delta > 0$ we denote by $\cl S_{\delta}$ the Riemann surface
obtained by the same construction as $\cl S_0$, pasting a plane
and cylinders isometrically along slits, but with ramification
points placed at $\dd Z - i\delta$ in the $0$-plane and $0 -
i\delta$ in the cylinders. The uniformization of this surface is
given by the map  $z \mapsto K(z+i\delta)$ instead of $z \mapsto
K(z)$. The image under this uniformization of the real line $\dd
R$ in the $0$-sheet of $\cl S_{\delta}$ is a 1-periodic graph with
an amplitude of oscillation that is large when $\delta$ is small.
For small $\delta$, the troughs of this curve have large negative
imaginary parts and lie well below the real axis; in applications
however, we would like the image of the real line to lie above the
$x$-axis, with the troughs touching the $x$-axis, and so we define
the uniformization $K_{\delta} : \cl S_{\delta} \to \dd C$ of $\cl
S_{\delta}$ by $$ K_{\delta}(z) := K(z+i\delta) - i \times
\hbox{"height of the troughs"} $$ where by the "height of the
troughs" we mean $$ \hbox{"height of the troughs"} \ = - \hbox{ Im }
K(n + i\delta) \ = {1 \over 2\pi } \log (-\log (1-e^{-2\pi
\delta}))$$ (where $n$ here is any integer). The curve
$K_{\delta}(\dd R)$ then lies in the upper half-plane of $\dd C$
with troughs lying on the real axis at the points $K_{\delta}(n +
0i ) = n + 0i, \, n \in \dd Z$; the image of the upper half-plane
of the $0$-sheet of $\cl S_{\delta}$ under $K_{\delta}$, which we
denote by $K_{\delta}(\dd H)$ is the region in the upper
half-plane of $\dd C$ bounded by $K_{\delta}(\dd R)$.

\medskip

In applications it will also be useful to have the image of the
real line oscillate at higher frequencies, more precisely with
period $1/a$ instead of period $1$ for integers $a \geq 1$; to do
this we consider the uniformizations ${1 \over a} K_{\delta}$
obtained by dividing $K_{\delta}$ by integers $a \geq 1$. The
image of the real line,  the curve ${1 \over a} K_{\delta}(\dd
R)$, is then a $1/a$-periodic curve. Finally, in order to control
the amplitude of oscillation, we note that given $a \geq 1$
integer and $h > 2$ real, there exists a unique $\delta =
\delta(a,h)$ (as is shown in [PM4] section 1.2.5), with $0 <
\delta < 1/2$, such that the amplitude of oscillation of
$K_{\delta}(\dd R)$ is $a \cdot h$, and hence that of ${1 \over a}
K_{\delta}(\dd R)$ is $h$. Moreover $\delta(a,h) \to 0$ for fixed
$h$ when $a \to \infty$. The constant $\delta$ satisfies (see
[PM4] section 1.2.5) $$ {\log \left(1 + e^{-2\pi \delta} \right)
\over -\log \left(1 - e^{-2\pi \delta} \right)} = e^{-2\pi ah} $$
The peaks and troughs of the curve ${1 \over a} K_{\delta}(\dd R)$
occur at the points ${1 \over a} K_{\delta}(n+1/2 + 0i) = n+1/2 +
ih$ and ${1 \over a} K_{\delta}(n + 0i) = n + 0i$ respectively,
and we have $\Max_{x \in \dd R} $ Im ${1 \over a} K_{\delta}(x) =
$ Im ${1 \over a} K_{\delta}(n+1/2) = h$, \ $\Min_{x \in \dd R} $
Im ${1 \over a} K_{\delta}(x) = $ Im ${1 \over a} K_{\delta}(n) =
0$, for any integer $n \in \dd Z$.

\medskip

{\hfill {\centerline {\psfig {figure=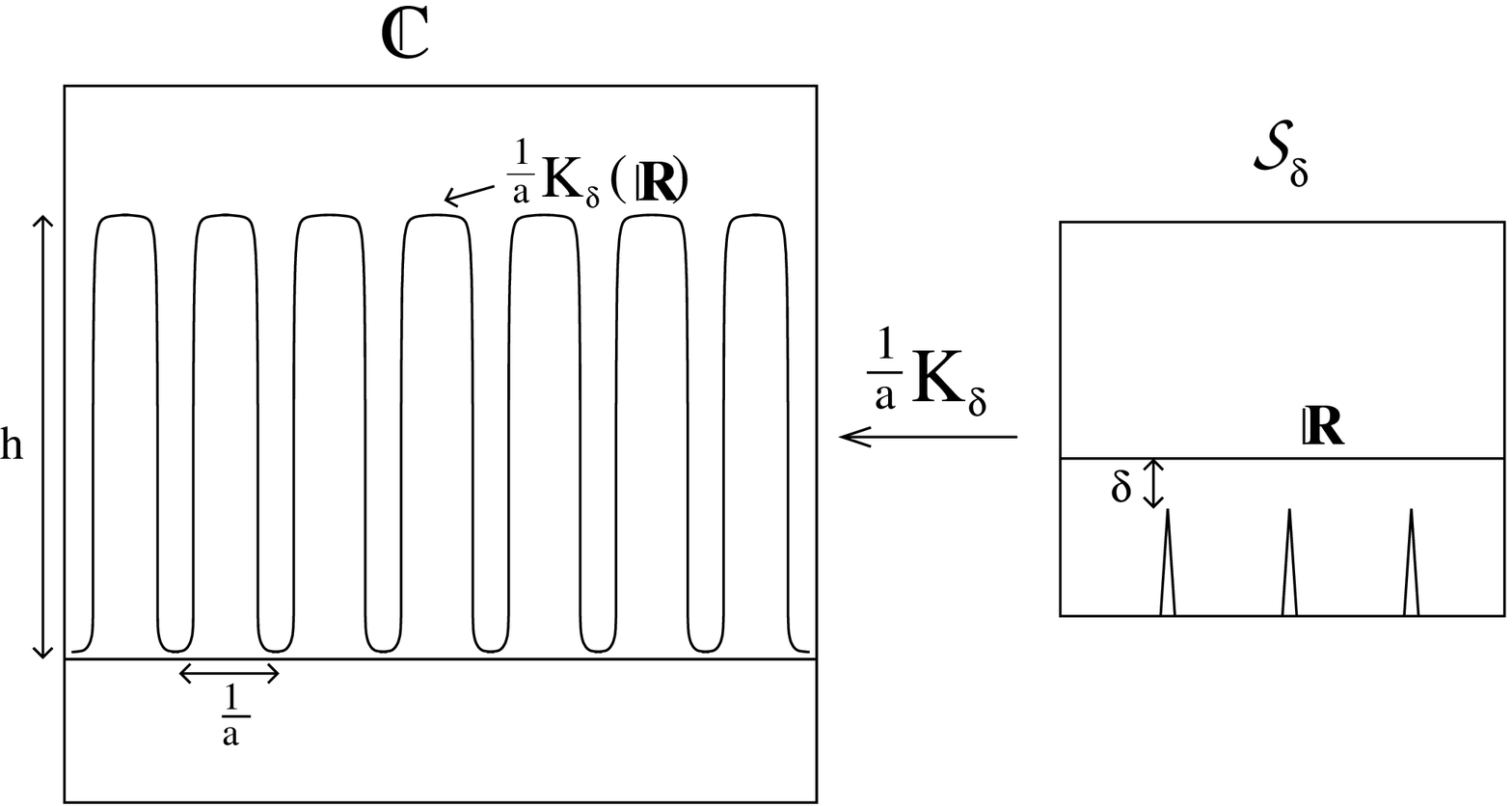,height=7cm}}}}

{\centerline {\bf Figure 3.}}

\medskip

We thus have a family of normalized uniformizations ${1 \over a}
K_{\delta}$ depending on two parameters $a \geq 1$ an integer and
$h > 2$ a real, where here $\delta$ is understood to be $\delta =
\delta(a,h)$. The maps ${1 \over a} K_{\delta}$ have the following
properties :

\medskip

{\bf (1)} ${1 \over a} K_{\delta}$ is univalent on a neighbourhood
of $\overline{ \dd{H} }$, mapping $\dd{H}$ into itself.

\medskip

{\bf (2)} ${1 \over a} K_{\delta}(z + 1) = {1 \over a}
K_{\delta}(z) + {1 \over a}$

\medskip

{\bf (3)} ${1 \over a} K_{\delta}(i \dd{R}_{\geq 0}) = i
\dd{R}_{\geq 0}, \, ({1 \over a} K_{\delta})'(z) > 0$ for $z \in i
\dd{R}_{\geq 0}$

\medskip

{\bf (4)} There exists a constant $C$ independent of $a$ and $h >
2$ such that $$ \left|{1 \over a} K_{\delta}(z) - \left({1 \over
a}z + ih \right) \right| \leq {C \over a} $$ for all $z$ with Im
$z > 1$.

\medskip

{\bf (5)} Let $z_0 \in \dd C$ denote the limit of ${1 \over a}
K_{\delta}(z)$ when $z \to \infty$ in $\cl S_{\delta}$ through the
upper end of one of the 1-cylinders. There exists a constant $C_1$
independent of $a$ and $h > 2$ such that $$ \hbox{Im } z_0 \geq h
- {C_1 \over a} $$ (note that the images under ${1 \over a}
K_{\delta}$ of the upper ends of the 1-cylinders differ from each
other by $1/a$-translations, so Im $z_0$ is independent of the
choice of 1-cylinder). We note that, for the dynamics obtained by
conjugation by ${1 \over a} K_{\delta}$, these points are no
longer necessarily fixed points (modulo $\dd{Z}$), but rather
periodic points (modulo $\dd{Z}$) of period $a$.

\medskip

We note that the properties {\bf (1)-(4)} refer to ${1 \over a}
K_{\delta}(z)$ for $z$ in the $0$-plane of $\cl{S}_{\delta}$,
while {\bf (5)} refers to the 1-cylinders. These properties appear
in [PM4], section 1.2.5, but can also be checked by elementary
computations . We note that condition {\bf (5)} allows us to
control the heights in the upper half-plane of the non-linearities
of the dynamics obtained by conjugation by ${1 \over a}
K_{\delta}$; indeed {\bf (5)} says that the heights of the
periodic orbits is approximately $h$, hence larger choices of $h$
result in somewhat smaller linearisation domains.

\medskip

{\bf 2.1.5)   Construction of the limit dynamics.}

\medskip

Now we iterate the procedure described in section 2.1.3 above. We
take a surface $\cl S_{\delta_1}$, lift the linear flow of
translations to $\cl S_{\delta_1}$ and conjugate by a
uniformization $K_0 = {1 \over a_0} K_{\delta_1}$ to obtain folded
dynamics as described in section 2.1.3. Having fixed $\cl
S_{\delta_1}$ and $K_0$, we take a surface $\cl S_{\delta_2}$,
lift the linear flow of translations to $\cl S_{\delta_2}$ and
conjugate by a uniformization $K_1 = {1 \over a_1} K_{\delta_2}$;
we then lift this folded dynamics to $\cl S_{\delta_1}$ and
conjugate by the uniformization $K_0$, to obtain a twice-folded
dynamics, given by conjugating by the map $K_0 K_1$. We continue
this way, choosing at each stage a surface $\cl S_{\delta_n}$ and
a uniformization $K_{n-1}$, to obtain at stage $n$ a dynamics
given by conjugation by $K_0 \dots K_{n-1}$. At each stage the
dynamics obtained is linearisable, but the linearisation domain
decreases; the point is that, by choosing the appropriate
uniformizations and controlling the geometry of the domains
obtained at each stage, we should be able to obtain, in the limit,
dynamics having an invariant domain with the desired geometry.

\medskip

In the construction that follows, we define a sequence $$K_n = {1
\over a_n} K_{\delta_{n+1}}, \ \delta_{n+1} = \delta(a_n, h_n), \
n \geq 0 $$ of uniformizations depending on two sequences of
parameters $(a_n)_{n \geq 0}$ and $(h_n)_{n \geq 0}$. The whole
construction is completely determined by the sequences of
parameters $a_n$ and $h_n$, which we will indicate how to choose
in the course of the construction. For the moment we assume that
they are given, and construct the desired limit dynamics as
follows:

\medskip

Let $n \geq 1$.

\medskip

We let $\cl{F}_{n,n} = (F_{n,n,t}(z) = z + t)_{t \in \dd{R}}$ be
the linear flow of translations. As explained in section 2.1.3
above, we can lift the elements of this flow to univalent maps
$\tilde{F}_{n,n,t}$ defined on the part of the surface
$\cl{S}_{\delta_{n}}$ obtained by deleting an
$\epsilon$-neighbourhood $\cl I(\epsilon)$ of the ramification
points and for times $t \in \dd Z + (-\epsilon, \epsilon)$. We
take $\epsilon = 3{a_n}^{-1}$, and conjugate by the uniformization
$K_{n-1}$ of $\cl S_{\delta_{n}}$ to obtain a semi-flow
$\cl{F}_{n-1,n} = (F_{n-1,n,t})_{t \in A_{n-1,n}}$ defined by $$
F_{n-1,n,{a_{n-1}}^{-1}t} := K_{n-1} \circ \tilde{F}_{n,n,t} \circ
K_{n-1}^{-1} $$ where $$ {a_{n-1}}^{-1}t \in A_{n-1,n} :=
{a_{n-1}}^{-1} (\dd Z + [-2{a_n}^{-1}, 2{a_n}^{-1}]). $$ (we could
in this step have taken a bigger set of times for $A_{n-1,n}$ such
as ${a_{n-1}}^{-1} (\dd Z + (-3{a_n}^{-1}, 3{a_n}^{-1}))$, but in
the steps that follow it will be necessary to restrict to smaller
sets of times). We note that the time has been reparametrized in
order to have the right behaviour at $+i\infty$, when Im $z \to
+\infty$, $$ F_{n-1,n,t}(z) = z + t + o(1). $$ These maps are
defined on $K_{n-1}(\cl{S}_{\delta_{n}} - \cl I(3{a_n}^{-1}))
\subset \dd C$, which contains a large half-plane $\dd H_{-y} = \{
\hbox{ Im } z \geq -y \}$, where $y \to +\infty$ as $a_n \to
+\infty$.

\medskip

We would now like to lift the semi-flow $\cl{F}_{n-1,n}$ to the
surface $\cl{S}_{\delta_{n-1}}$ and conjugate by the
uniformization $K_{n-2}$ to obtain a semi-flow $\cl{F}_{n-2,n}$.
The difference with the previous case is that we are now starting
with a nonlinear flow instead of the linear flow of translations.
However, by giving up some space in the lower region of the
half-plane of definition $\dd H_{-y}$, we can work with maps close
to linear ones; more precisely we have (see [PM4] Lemma 1.1.6 and
Proposition 2.2.6), for $t \in A_{n-1,n}$, if Im $z \geq -y + {1
\over 2\pi} \log a_{n-1} + C$, (where $C$ is a universal constant)
then $$ |F_{n-1,n,t}(z) - (z + t)| < {1 \over a_{n-1}} $$ Since $y
\to +\infty$ as $a_n \to +\infty$, we observe that once $a_{n-1}$
is fixed, by choosing $a_n$ large enough we can ensure that this
inequality holds in as large a half-plane $\dd H_{-\tilde{y}}$ as
desired, where $-\tilde{y} = -y + {1 \over 2\pi} \log a_{n-1} +
C$.

\medskip

Now consider the region $\cl J(-\tilde{y}) = \{ z \in \cl
S_{\delta_{n-1}} : \hbox{ Im }z \leq -\tilde{y} \} \subset \cl
S_{\delta_{n-1}}$ (note that there is a well-defined function Im
$z$ on $\cl S_{\delta_{n-1}}$). Since the estimate above holds on
the complement $\cl S_{\delta_{n-1}} - \cl J(-\tilde{y})$ of this
region, it should be clear that if we now remove as well an
$\epsilon$-neighbourhood $\cl I(\epsilon)$ of the ramification
points, this time with $\epsilon = 3{a_{n-1}}^{-1}$, then we can
define, for $t \in \dd  Z + A_{n-1,n} \cap [-2{a_{n-1}}^{-1},
2{a_{n-1}}^{-1}]$, the lifts of elements of $\cl{F}_{n-1,n}$ as
univalent maps $\tilde{F}_{n-1,n,t} : \cl S_{\delta_{n-1}} - (\cl
J(-\tilde{y}) \cup \cl I(3{a_{n-1}}^{-1})) \to \cl
S_{\delta_{n-1}}$.

\medskip

We can now proceed as before and define the semi-flow
$\cl{F}_{n-2,n} = (F_{n-2,n,t})_{t \in A_{n-2,n}}$ by $$
F_{n-2,n,{a_{n-2}}^{-1}t} := K_{n-2} \circ \tilde{F}_{n-1,n,t}
\circ K_{n-2}^{-1} $$ where $$ {a_{n-2}}^{-1}t \in A_{n-2,n} :=
{a_{n-2}}^{-1} (\dd Z + A_{n-1,n} \cap [-2{a_{n-1}}^{-1},
2{a_{n-1}}^{-1}]) $$ (the times have again been reparametrized
appropriately so that $F_{n-2,n,t}(z) = z+t+o(1)$ when Im $z \to
+\infty$).

\medskip

These maps are defined on $K_{n-2}\left(\cl S_{\delta_{n-1}} -
(\cl J(-\tilde{y}) \cup \cl I(3{a_{n-1}}^{-1}))\right) \subset \dd
C$, which, as remarked earlier in section 2.1.3, contains a large
half-plane $\dd H_{-y'} = \{ \hbox{ Im } z \geq -y' \}$, where $y'
\to +\infty$ as $a_{n-1} \to +\infty$ and $\tilde{y} \to \infty$
(see [PM4] Proposition 2.2.8).

\medskip

Following this same procedure of lifting and conjugating the
dynamics, we can define successively the semi-flows
$\cl{F}_{n-3,n}, \, \cl{F}_{n-4,n}, \dots, \, \cl{F}_{0,n}$, with
$\cl{F}_{j,n} = (F_{j,n,t})_{t \in A_{j,n}}$ defined for times $t
\in A_{j,n}$ where $A_{j,n}$ is defined inductively by $$ A_{j,n}
:= {a_{j}}^{-1} (\dd Z + A_{j+1,n} \cap [-2{a_{j+1}}^{-1},
2{a_{j+1}}^{-1}]) $$ for $j=n-1,\dots,0$ and $A_{n,n} = \dd R$. We
note here that the sets $A_{0,n}$ form a decreasing sequence, ie
$A_{0,n+1} \subset A_{0,n}, n \geq 1$; indeed, $$ A_{n-1,n+1} =
{a_{n-1}}^{-1} (\dd Z + A_{n,n+1} \cap [-2{a_n}^{-1},
2{a_n}^{-1}]) \subset {a_{n-1}}^{-1} (\dd Z + [-2{a_n}^{-1},
2{a_n}^{-1}]) = A_{n-1,n} $$ from which it follows similarly that
$$\eqalign{ A_{n-2,n+1} & \subset A_{n-2,n} \cr A_{n-3, n+1} &
\subset A_{n-3,n} \cr & \dots \cr A_{0,n+1} & \subset A_{0,n} \cr
}$$.

Thus we obtain a sequence of semi-flows $\cl{F}_{0,n} =
(F_{0,n,t})_{t \in A_{0,n}}$, where each one is given by 'folding'
the linear flow $n$ times. Moreover, the $a_n$'s can be chosen
inductively to grow fast enough so that all the semi-flows
$\cl{F}_{0,n}$ are defined in the upper half-plane $\dd H$, or
indeed in any fixed half-plane $\dd H_{-M} = \{ \hbox{ Im } z > -M
\}$, as large as desired; indeed we can state as a proposition the
following (which follows directly from [PM4] sections 1.3.1 and
2.2.3):

\medskip

\eno{Proposition 2.1.1}{ Given $M > 0$, there exists a sequence of
conditions $(\cl C_n)_{n \geq 1}$ on the integers $a_n$ of the
form $$ (\cl C_n) \qquad \quad a_n \geq C_n(a_0, \dots, a_{n-1},
h_0, \dots, h_{n-1}, M) $$ such that when they are fulfilled then
all the semi-flows $(\cl{F}_{0,n})_{n \geq 1}$ are defined in the
half-plane $\dd H_{-M} = \{ \hbox{ Im } z > -M \}$.}

\medskip

We observe thus a feature of the construction which is crucial
when passing to the limit, namely that it allows us to
successively paste nonlinearities to the dynamics without
decreasing its domain of definition too much, allowing us to
retain a half-plane $\dd H_{-M}$.

\medskip

Assume that the $a_n$'s have been chosen so that all the
semi-flows $\cl{F}_{0,n}$ are defined on $\dd H_{-M}$. From the
sequence of semi-flows $\cl{F}_{0,n}$ we can extract a convergent
subsequence (it always exists, see [PM4], section 1.1.2)
converging normally for times $t \in A = \cap_{n \geq 0} A_{0,n}$
to get a limit semi-flow $\cl{F} = (F_t)_{t \in A}$ defined on the
half-plane $\dd H_{-M}$ for times $t \in A$. We note that elements
of the limit semi-flow do indeed commute, since we can pass to the
limit in the equation $$F_{0,n_k,t} \circ F_{0,n_k,s} (z) =
F_{0,n_k,s} \circ F_{0,n_k,t} (z) $$ (using uniform convergence in
compact neighbourhoods of the points $z, F_t(z), F_s(z)$).

\medskip

The Cantor set of times $C$ in the main theorem is given by
$$\displaystyle{ C := \bigcap_{n \geq 0} A_{0,n} / \dd Z }$$ If
all the $a_n$'s are large enough, say for example $a_n \geq 5$ for
all $n$, then the set $C$ is indeed a Cantor set.

\medskip

The semi-flow $(f_t)_{t \in C}$ of the main theorem is defined on
a common neighbourhood $E(\dd H_{-M})$ of the closed unit disk
$\overline{\dd{D}}$ via the equation $$ f_t \circ E = E \circ F_t,
\ t \in C $$ (where $E$ is the universal covering $E : \dd C \to
\dd C - \{0\}, E(z) = e^{2\pi i z}$).

\medskip

As stated in section 2.1.4, each uniformization $K_n$ maps the
upper half-plane of the $0$-plane of $\cl S_{\delta_{n+1}}$ into
the upper half-plane $\dd H$. So we can consider its restriction
to the upper-half plane of the $0$-plane as a map $K_n : \dd H \to
\dd H$ of $\dd H$ into itself. If we consider the composition of
these maps $K_0 \circ K_1 \circ \dots \circ K_n$, then it is clear
that the set $K_0 K_1 \dots K_n (\overline{\dd{H}})$ is totally
invariant under the semi-flow $\cl{F}_{0,n+1}$, and hence the set
$$ \cl{H} = \bigcap_{n = 0}^{\infty} K_0 K_1 \dots K_n
(\overline{\dd{H}}) $$ is then totally invariant under the limit
semi-flow $\cl{F}$. Moreover, $\cl{H}$, being the decreasing limit
of the domains $K_0 K_1 \dots K_n (\overline{\dd{H}})$, is full,
and $\cl{H} \cap \partial \dd{H} \neq \phi$ (since $0 \in
\cl{H}$). Thus $\cl{H}$ is in fact a common hedgehog for
irrational elements of the limit semi-flow $\cl{F}$. The hedgehog
$K$ of the semi-flow $(f_t)_{t \in C}$ in the main theorem is
given by $$ K = E(\cl{H}) $$

\medskip

It is important to note in this construction that the only
requirement imposed on the $a_n$'s so that the limit dynamics is
defined on the upper half-plane is that they grow fast enough, as
expressed by Proposition 2.1.1 above. Subject to this restriction,
the construction still offers enough flexibility to construct
hedgehogs with quite different geometries by varying the choices
of the parameters $a_n$ and $h_n$.

\medskip

In the following sections we show how to choose the
uniformizations $K_n$ appropriately (by proper choices of the
parameters $a_n, h_n$), and control the geometries of the domains
$K_0 K_1 \dots K_n(\overline{\dd{H}})$, so that the hedgehog
$\cl{H}$ has Hausdorff dimension 1.

\bigskip

{\bf 2.2  A hedgehog of Hausdorff dimension 1.}

\medskip

We fix some notation first; for a countable family $\cl{A} = \{ U_i \}_{i \in I}$ of planar sets $U_i$, for
$s > 0$ its $s$-$length$ is defined to be
$$
\cl{L}^s(\cl{A}) = \sum_{i \in I} \hbox{diam}(U_i)^s
$$
(recall that this is the quantity used to define Hausdorff $s$-measure).

\medskip

We recall the definition of the uniformization ${1 \over a}
K_{\delta}$,
$$
{1 \over a}K_{\delta}(z) = {1 \over a} \left( {1 \over 2 \pi i} \log \left( - \log( 1 - e^{2 \pi i (z + i\delta)} )
\right)- iC_{\delta} \right)
$$
where $C_{\delta} = {1 \over 2\pi } \log (-\log (1-e^{-2\pi
\delta}))$. We consider the part of the domain ${1 \over a} K_{\delta}
(\overline{\dd{H}})$ below the height $h - 1$.

\medskip

{\hfill {\centerline {\psfig {figure=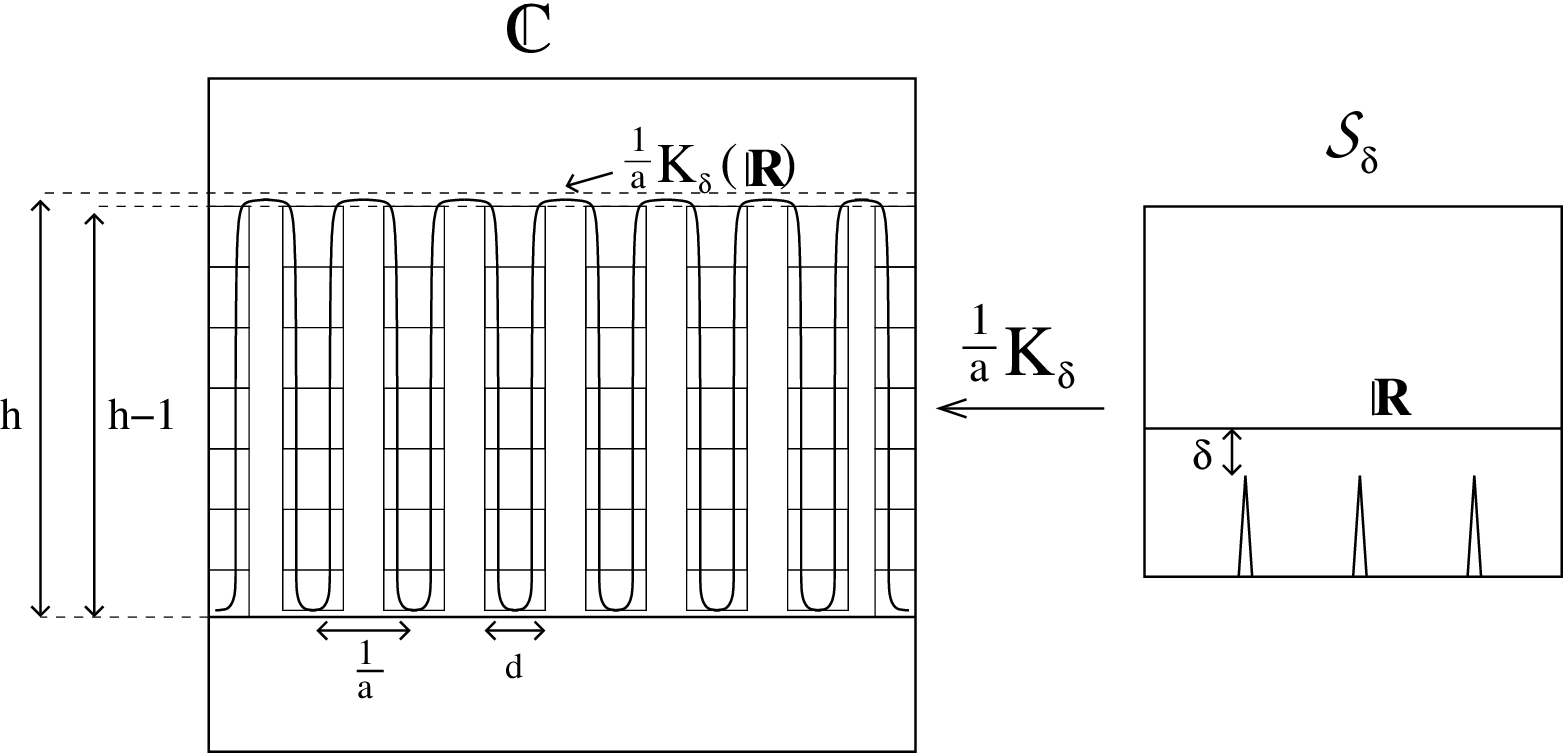,height=7cm}}}}

{\centerline {\bf Figure 4.}}

\medskip

For the parametrization $t \mapsto {1 \over a} K_{\delta}(t) = x(t) + iy(t)$ of the curve ${1 \over a}
K_{\delta}(\dd{R})$, it is straightforward to check from the formulae that

\medskip

(1) $x'(t) > 0$ for all $t$, so the curve is a graph over $\dd{R}$, ie has a parametrization of the form $y =
f(x)$.

\medskip

(2) Since ${1 \over a} K_{\delta}(t+1) = {1 \over a} K_{\delta}(t) + {1 \over a}$, the graph is ${1 \over a}$
periodic, ie $f(x+{1 \over a}) = f(x)$. Moreover in any period $n \leq t \leq n+1, \, n \in \dd{Z}$, we have
$$\displaylines{
x(n) + iy(n) = {n \over a} + i0,  \quad x(n+1/2) + iy(n+1/2) = {n+1/2 \over a} + ih,  \cr
y'(t) > 0 \hbox{ for } n < t < n+1/2, \quad y'(t) < 0 \hbox{ for } n+1/2 < t < n+1, \hbox{ and}  \cr
y((n+1/2)-u)=y((n+1/2)+u) \quad \hbox{ for all } u \in \dd{R} \cr
}$$

So in every period $n \leq t \leq n+1$, there are exactly two points $t = n+t_0$ and $t = n+1-t_0, \, 0 < t_0
< 1$, where the curve crosses the line $\{ $ Im $z = h - 1 \}$. Moreover $0 \leq y(t) \leq h-1$ for $n-t_0
\leq t \leq n+t_0$. Thus taking $d = 2 \, x(t_0)$, we can, given an integer $M > 0$, cover the region ${1
\over a} K_{\delta}(\overline{\dd{H}}) \cap \{ \hbox{ Im } z \leq h-1 , \, 0 \leq \hbox{ Re } z \leq M \} =
\{ z=x+iy : f(x) \leq y \leq h-1, \, 0 \leq x \leq M \}$ by squares of sidelength $d$, as shown in figure 4
above.

\medskip

We denote this covering by $\cl{A}(a,h,M)$.

\medskip

The key to the whole construction is the following lemma, which says that the region being covered can be
made arbitrarily 'thin' by taking $a$ large enough, keeping $h,M$ fixed:

\medskip

\eno {Lemma.} {Fix $s>1$. With the same notation as above, we have
$$
d = O \left( {1 \over a} e^{-2 \pi a} \right) \, \, \hbox{ and } \, \, \cl{L}^s(\cl{A}(a,h,M)) \to 0 \, \, \,
\hbox{ as } a \to +\infty, \hbox{ keeping } h,M \hbox{ fixed.}
$$
}

\medskip

\stit {Proof.} We first prove the estimate on $d = 2 \, x(t_0)$.

$t_0$ was defined as the smallest $t>0$ such that $y(t)=$ Im ${1
\over a} K_{\delta}(t) = h-1$. The formula for ${1 \over a}
K_{\delta}$ gives $$ {1 \over 2 \pi} \log \left| { \log(1 - e^{-2
\pi \delta}) \over \log(1 - e^{2 \pi i (t_0 + i\delta)}) } \right|
= a(h-1) $$ But $\delta = \delta(a,h)$ is chosen such that $$ {1
\over 2 \pi} \log \left| { \log(1 - e^{-2 \pi \delta}) \over
\log(1 + e^{-2 \pi \delta)}) } \right| = ah \, ; $$ substituting
this into the previous equation, rearranging terms and
exponentiating gives $$ | \log(1 - e^{2 \pi i (t_0 + i\delta)}) |
= e^{2\pi a} \log(1 + e^{-2\pi \delta}) $$ For fixed $h > 0$, as
$a \to +\infty$, $\delta \to 0$ (see [PM4]), so $$ | \log(1 - e^{2
\pi i (t_0 + i\delta)}) | \sim e^{2\pi a} \log(2) $$ as $a \to
+\infty$. Thus $$\eqalign{ d & = 2 \, x(t_0) = 2 \hbox{ Re } {1
\over a} K_{\delta}(t_0) \cr
  & = 2 \, {1 \over a} {1 \over 2 \pi} \arg \left( - \log( 1 - e^{2 \pi i (t_0 + i\delta)} ) \right) \cr
  & = 2 \, {1 \over a} {1 \over 2 \pi} \arcsin \left( { \hbox{Im }w \over |w| } \right) \qquad (\hbox{where }
w = - \log( 1 - e^{2 \pi i (t_0 + i\delta)} )  \cr
  & = 2 \, {1 \over a} {1 \over 2 \pi} O \left( { \hbox{Im }w \over |w| } \right) \qquad \qquad (\hbox{note
that } \left( { \hbox{Im }w \over |w| } \right) \to 0 \, \hbox{ in this case}) \cr
  & = O \left( {1 \over a} e^{-2 \pi a} \right) \qquad \qquad \quad (\hbox{since Im }w \hbox{ stays bounded,
and } |w| \sim e^{2\pi a} \log(2) ) \cr
}$$
The covering $\cl{A}(a,h,M)$ has a number $N$ of squares of the order of $M \cdot {h \over d} \cdot a$, so
$$\eqalign{
\cl{L}^s(\cl{A}(a,h,M)) & = N \cdot (\sqrt{2} d)^s \cr
                      & \sim M \cdot {h \over d} \cdot a \cdot (\sqrt{2} d)^s \cr
                      & = M \cdot h \cdot O \left( \left({1 \over a} e^{-2 \pi a} \right)^{s-1} \cdot a
\right) \cr
                      & \to 0 \qquad \hbox{ as } a \to +\infty \qquad \qquad \qquad \qquad \qquad \qquad
\qquad \qquad \qquad \diamondsuit. \cr
}$$

\medskip

{\it Remark.} To conclude $\cl{L}^s(\cl{A}(a,h,M)) \to 0$ as $a
\to +\infty$, one only needs in fact $d = o(1/a)$ and not the
stronger exponential decay proved above.

\medskip

The idea of the proof is now simple: we choose the maps $K_n$ inductively. Say $K_0, \dots K_{n-1}$ have been
chosen; by fixing $h_n$ to be quite large, and then choosing $a_n$ very large, the lemma above ensures that a
large part of the domain $K_n(\overline{\dd{H}})$ can be made very 'thin'. Since $K_0 \dots K_{n-1}$ is a
fixed holomorphic diffeomorphism not depending on the choice of $a_n$, we can ensure by choosing $a_n$ large
enough, that it stays thin even after transporting by $K_0 \dots K_{n-1}$. Thus at each stage $n$, a large
part of the domain $K_0 \dots K_n(\overline{\dd{H}})$ can be made very thin, so making the domains thinner
and thinner at each stage we can obtain in the limit a hedgehog with Hausdorff dimension 1.

\medskip

\stit{Proof of Main Theorem.}

Fix three decreasing sequences $s_n \downarrow 1, \  \epsilon_{n} \downarrow 0$ and $d_n \downarrow 0$.
We choose the parameters $a_n$ and $h_n, n \geq 0$, inductively as follows:

\medskip

For $n = 0$ : \quad  Let $h_0 = 10$. Applying the above lemma to $K_0 = {1 \over a_0} K_{\delta(a_0, h_0)}$,
we choose $a_0 \geq 1$ large enough so that
$$\eqalign{
\hbox{diam}(U) & \leq d_0 \qquad \forall\ U \in \cl{A}(a_0,h_0,1), \, \hbox{ and } \cr
\cl{L}^{s_0}(\cl{A}(a_0,h_0,1)) & \leq \epsilon_0. \cr
}$$

For $n \geq 1$: \quad  Assume $a_0, \dots, a_{n-1}$ and $h_0, \dots, h_{n-1}$ have been chosen, so the
functions $K_0, \dots, K_{n-1}$ are determined and fixed.

\medskip

By property {\bf (4)} of the maps ${1 \over a} K_{\delta}$, we
have $K_0 \circ \dots \circ K_{n-1}(z) = z + i(h_0 + \dots +
h_{n-1}) + O(1/a_0 \dots a_{n-1})$ as Im $z \to +\infty$. We
choose $h_n > 0$ large enough so that Im $K_0 \dots K_{n-1}(x +
i(h_n - 1)) \geq h_0 + n$ for all $x$.

\medskip

By the above choice of $h_n$, $K_0 \dots K_{n-1} (\overline{\dd{H}}) \cap \{ \hbox{ Im } z \leq h_0 + n \}$
is contained in $K_0 \dots K_{n-1} (\{ 0 \leq \hbox{ Im } z \leq h_n - 1 \})$. Thus the covering $\cl{A}(a_n,
h_n, \, a_0 \dots a_{n-1})$ of the region $K_n(\overline{\dd{H}}) \cap \{ \hbox{ Im } z \leq h_n - 1 , \, 0
\leq \hbox{ Re } z \leq a_0 \dots a_{n-1} \}$ is transported by $K_0 \dots K_{n-1}$ to a covering
$$
\cl{A}_n = K_0 \dots K_{n-1} \left( \cl{A}(a_n,h_n, \, a_0 \dots a_{n-1}) \right)
$$
of the region $K_0 \dots K_n (\overline{\dd{H}}) \cap \{ \hbox{ Im } z \leq h_0 + n , \, 0 \leq \hbox{ Re } z
\leq 1 \}$. By the lemma, for the covering $\cl{A}(a_n,h_n, \, a_0 \dots a_{n-1})$, for fixed $s$ both its
$s$-$length$ and the maximum diameter of the sets it contains go to $0$ as $a_n \to +\infty$. Since $K_0
\dots K_{n-1}$ is a fixed holomorphic diffeomorphism independent of the choice of $a_n$, the same is true for
the covering $\cl{A}_n$. We thus choose $a_n$ large enough so that

$$\eqalign{
\hbox{diam}(U) & \leq d_n \qquad \forall\ U \in \cl{A}_n, \, \hbox{ and } \cr
\cl{L}^{s_n}(\cl{A}_n) & \leq \epsilon_n. \cr
}$$

Finally, we also choose $a_n$ large enough depending on
$a_0,\dots,a_{n-1},h_0,\dots,h_n$ and as given by Proposition
2.1.1 so that we do obtain limit dynamics defined on a half-plane
strictly containing $\overline{\dd{H}}$.

\medskip

We can now prove that with the choice of the parameters $a_n, \, h_n, \, n \geq 0$ as above, the resulting
hedgehog
$$
\cl{H} = \bigcap_{n = 0}^{\infty} K_0 \dots K_n(\overline{\dd{H}})
$$
has Hausdorff dimension 1 :

\medskip

Its enough to show that $H^s(\cl{H}) = 0$ for all $s > 1, \, H^s$ denoting $s$-dimensional
Hausdorff measure. So fix $s>1$.

\medskip

For $M > 0$, writing $\cl{H}_M = \cl{H} \cap \{ \hbox{ Im } z \leq M, \, 0 \leq \hbox{ Re } z \leq 1 \}$, we
have
$$\eqalign{
H^s(\cl{H}_M) & = \lim_{\delta \to +0} \ H_{\delta}^s(\cl{H}_M) \cr
 & = \lim_{n \to \infty} \ H_{d_n}^s(\cl{H}_M) \qquad \qquad \quad (\hbox{ since } d_n \to 0 ) \cr
 & \leq \limsup_{n \to \infty} \ H_{d_n}^{s_n}(\cl{H}_M) \qquad \qquad (\hbox{ since } s_n \leq s \hbox{
eventually, } H_{\delta}^s \hbox{ decreasing in } s) \cr
 & \leq \limsup_{n \to \infty} \ \cl{L}^{s_n}(\cl{A}_n) \qquad \qquad (\hbox{ since } \cl{A}_n \hbox{ covers
} \cl{H}_M \hbox{ once } h_0 + n \geq M ) \cr
 & \leq \limsup_{n \to \infty} \ \epsilon_n \cr
 & = 0 \cr
}$$
So
$$\eqalign{
H^s(\cl{H} \cap \{ 0 \leq \hbox{ Re } z \leq 1 \}) & = \lim_{M \to \infty} \ H^s(\cl{H}_M) \cr
 & = 0 \cr
}$$
Since $\cl{H} = \cl{H} + 1$, it follows that $H^s(\cl{H}) = 0$, and hence $H^s(K) = 0$. \qquad \qquad \qquad
\qquad $\diamondsuit$.

\bigskip

\stit {3. Remarks.}

\bigskip

1. Note that since the hedgehog has Hausdorff dimension 1, it has no interior and cannot contain a Siegel
disk, so the limit dynamics $(f_t)_{t \in C}$ is nonlinearisable.

\medskip

2. The rotation numbers $t \in C$ of the limit dynamics are of the form
$$
t = \sum_{n = 0}^{\infty} {{\varepsilon}_n \over {a_0 \dots a_n}}
$$
where each ${\varepsilon}_n = 0, \, 1$ or $-1$. Since our construction above requires the $a_n$'s to be
growing very fast, this means that we obtain very Liouville rotation numbers.

\medskip

3. It is possible by modifying the construction appropriately to obtain linearisable dynamics having a common
Siegel disk, but such that the portion of the hedgehog outside the closure of the Siegel disk has Hausdorff
dimension 1. Thus linearisability is not an obstruction to the hedgehog being 'thin'.

\bigskip

{\bf 4. Acknowledgements.} The present article is based on a
portion of my PhD thesis. I thank my advisor Ricardo Perez-Marco
for his guidance and for introducing me to his techniques of
tube-log Riemann surfaces.

\bigskip

{\centerline {\bf References}}

\bigskip

\medskip

[Bi] G.D.BIRKHOFF, {\it Surface transformations and their dynamical applications}, Acta Mathematica, {\bf 43},
(1920). Also, Collected Mathematical Papers II, p.111

\medskip

[B] K.BISWAS, {\it Smooth combs inside hedgehogs}, Discrete and
Continuous Dynamical Systems, vol 12, no.5, May 2005, p.853-880.

\medskip

[Cr1] H.CREMER, {\it Zum Zentrumproblem}, Math Ann., vol 98, 1928, p. 151-153

\medskip

[Cr2] H.CREMER, {\it \"Uber die H\"aufigkeit der Nichtzentren}, Math Ann., vol 115, 1938, p.573-580

\medskip

[PM1] R.PEREZ-MARCO, {\it Fixed points and circle maps}, Acta Mathematica, {\bf 179:2}, 1997, p.243-294.

\medskip

[PM2] R.PEREZ-MARCO,  {\it Sur les dynamiques holomorphes non
lin\'earisables et une conjecture de V.I. Arnold}, Ann. Scient.
Ec. Norm. Sup. 4 serie, {\bf 26}, 1993, p.565-644; (C.R. Acad.
Sci. Paris, {\bf 312}, 1991, p.105-121).

\medskip

[PM3] R.PEREZ-MARCO, {\it Uncountable number of symmetries for
non-linearisable holomorphic dynamics}, Inventiones Mathematicae,
{\bf 119}, {\bf 1}, p.67-127,1995; (C.R. Acad. Sci. Paris, {\bf
313}, 1991, p. 461-464).

\medskip

[PM4] R.PEREZ-MARCO, {\it Siegel disks with smooth boundary},
Preprint, Universit\'e de PARIS-SUD, 1997,
www.math.ucla.edu/$\tilde{\ }$ricardo, submitted to the Annals of
Mathematics in 1997.

\medskip

[PM5] R.PEREZ-MARCO, {\it Hedgehog's Dynamics}, Preprint, UCLA, 1996.

\medskip

[PM6] R.PEREZ-MARCO,  {\it Sur les dynamiques holomorphes non lin\'earisables et une conjecture de V.I. Arnold},
Ann. Scient. Ec. Norm. Sup. 4 serie, {\bf 26}, 1993, p.565-644; (C.R. Acad. Sci. Paris, {\bf 312}, 1991,
p.105-121).

\medskip

[Si] C.L.SIEGEL, {\it Iteration of Analytic Functions}, Ann. Math., vol.43, 1942, p. 807-812

\medskip

[Y] J.C.YOCCOZ, {\it Petits diviseurs en dimension 1}, S.M.F., Asterisque {\bf 231} (1995)

\end